\documentclass[11pt]{article}
\usepackage[dvips]{graphics}
\usepackage{epsfig,rotating}
\usepackage{amssymb,amsmath}
\usepackage{latexsym}
\usepackage[usenames]{color}

\newcommand{\be}{\begin{equation}}
\newcommand{\ee}{\end{equation}}
\newcommand{\beaa}{\begin{eqnarray*}}
\newcommand{\eeaa}{\end{eqnarray*}}
\newcommand{\bea}{\begin{eqnarray}}
\newcommand{\eea}{\end{eqnarray}}
\newcommand{\lbl}{\label}

\newcommand{\bei}{\begin{itemize}}
\newcommand{\eei}{\end{itemize}}

\newcommand{\bd}{\bold}

\newtheorem{theorem}{ \noindent T{\footnotesize HEOREM}}

\newtheorem{lemma}{ \noindent L{\footnotesize EMMA}}[section]
\newtheorem{coro}{ \noindent C{\footnotesize OROLLARY}}[section]
\newtheorem{remark}{ \noindent R{\footnotesize EMARK}}[section]

\newcommand{\goto}{\rightarrow}
\begin{document}

\title{Phase Transition in Limiting Distributions of Coherence of High-Dimensional Random Matrices}
\author{Tony Cai$^{1}$ and Tiefeng Jiang$^{2}$\\
University of Pennsylvania and University of Minnesota}

\date{}
\maketitle

\footnotetext[1]{ Statistics Department, The Wharton School, University of
  Pennsylvania,  Philadelphia, PA 19104, \newline  \indent \ \
tcai@wharton.upenn.edu. The research of Tony Cai was supported in part
by NSF FRG Grant  \newline  \indent \ \
DMS-0854973.}
\footnotetext[2]{School of Statistics, University of Minnesota, 224 Church
Street, MN55455, tjiang@stat.umn.edu. \newline  \indent \ \
The research of Tiefeng Jiang was
supported in part by NSF FRG Grant DMS-0449365.}

\begin{abstract}
The coherence of a random matrix, which is  defined to be the largest magnitude of the Pearson correlation coefficients between the columns of  the random matrix, is an important quantity for a wide range of  applications including high-dimensional statistics and signal processing. Inspired by these applications, this paper studies the limiting laws of the coherence of $n\times p$ random matrices for a full range of the dimension $p$ with a special focus on the ultra high-dimensional setting. Assuming the  columns of the random matrix are independent random vectors with a common spherical distribution, we give a complete characterization of the behavior of the limiting distributions  of  the coherence. More specifically,  the limiting distributions of the coherence are derived separately for three regimes: $\frac{1}{n}\log p \to 0$, $\frac{1}{n}\log p \to \beta\in (0, \infty)$, and $\frac{1}{n}\log p \to\infty$.  The results show that the limiting behavior of the coherence differs significantly in different regimes and exhibits interesting phase transition phenomena as the dimension $p$ grows as a function of $n$. Applications to statistics  and compressed sensing in the ultra high-dimensional setting are also discussed.
\end{abstract}




\newpage
\section{Introduction}
\lbl{intro}
\setcounter{equation}{0}

With dramatic advances in computing and technology, large and high-dimensional datasets are now routinely collected in many scientific investigations. The associated statistical inference problems, where the dimension $p$ can be much larger than the sample size $n$, arise naturally in a wide range of applications including compressed sensing, climate studies, genomics, functional magnetic resonance imaging, risk management and portfolio allocation. Conventional statistical methods and results based on fixed $p$ and large $n$ are no longer applicable and these applications call for new technical tools and new statistical procedures.

The coherence of a random matrix, which is defined to be the largest magnitude of the off-diagonal entries of the sample correlation matrix generated from the random matrix, has been shown to be an important quantity for many applications. For example, the coherence has been used for testing the covariance structure of high-dimensional distributions (Cai and Jiang (2010)), the construction of compressed sensing matrices  and high dimensional regression in statistics (see, e.g.,  Candes and Tao (2005), Donoho, Elad and Temlyakov (2006) and Cai, Wang and Xu (2010a, b)). In addition, the coherence has also been used in signal processing, medical imaging, and seismology. Some of these problems are seemingly unrelated at first sight, but interestingly they can all be attacked through the use of the limiting laws of the coherence of random matrices (see, e.g., Cai and Jiang (2010)). In these applications, a case of special interest is when the dimension $p$ is much larger than the sample size $n$. Indeed, in compressed sensing and other related problems the goal is often to make the dimension $p$ as large as possible relative to the sample size $n$.

In the present paper  we study the limiting laws of the coherence of random matrices. Let $\bd{x}=(x_1, \cdots, x_n)^T\in \mathbb{R}^n$ and  $\bd{y}=(y_1, \cdots, y_n)^T\in \mathbb{R}^n.$ Recall the Pearson correlation coefficient $\rho$ defined by
\bea\lbl{tube}
\rho=\rho_{\bd{x},\bd{y}}=\frac{\sum_{i=1}^n(x_i-\bar{x})(y_i - \bar{y})}
{\sqrt{\sum_{i=1}^n(x_i-\bar{x})^2\cdot \sum_{i=1}^n(y_i-\bar{y})^2}}
\eea
where $\bar{{x}}=\frac{1}{n}\sum_{i=1}^nx_i$ and $\bar{{y}}=\frac{1}{n}\sum_{i=1}^ny_i$. Let $\bd{X}_1, \cdots, \bd{X}_p$ be independent $n$-dimensional random vectors, and let $\rho_{ij}$ be the correlation coefficient between $\bd{X}_i$ and $\bd{X}_j.$ Set $\bd{X}=(\bd{X}_1, \cdots, \bd{X}_p)=(x_{ij})_{n\times p}$. The coherence of the random matrix  $\bd{X}$ is defined as
\bea\lbl{mount}
L_n=\max_{1\leq i <  j \leq p}|\rho_{ij}|.
\eea
In certain applications such as the construction of compressed sensing
matrices, the means $\mu_i = E \bd{X}_i$ and $\mu_j = E\bd{X}_j$ are given and one is  interested in
\be\lbl{corr'}
\tilde{\rho}_{ij}=\frac{(\bd{X}_i-\mu_i)^T(\bd{X}_j-\mu_j)}{\|\bd{X}_i-\mu_i\|\cdot \|\bd{X}_j-\mu_j\|},\
\ 1\leq i, j\leq p
\ee
and the corresponding coherence is defined by
\bea\lbl{mount3}
\tilde{L}_n=\max_{1\leq i <  j \leq p}|\tilde{\rho}_{ij}|.
\eea
The goal of this paper is to give a complete characterization of the behavior of the limiting distributions  of $L_n$ and $\tilde{L}_n$ over the full range of $p$ (as a function of $n$) including the super-exponential case where $(\log p)/n \to \infty$.

The coherence $L_n$ has been studied intensively in recent years. Jiang (2004) was the first to show that if $x_{ij}$'s are independent and identically distributed (i.i.d.) with $E|x_{ij}|^{30+\epsilon} <\infty$ for some $\epsilon >0$ and $n/p \to \gamma\in (0, \infty),$   then $nL_n^2-4\log p + \log\log p$ converges weakly to an extreme distribution of type I with distribution function
\be\lbl{preceding}
F(y)=e^{-{1\over \sqrt{8\pi}}e^{-y/2}}, \;\; y\in \mathbb{R}.
\ee
Throughout this paper, $\log x=\log_e x$ for any $x>0$ and $p=p_n$ depends on $n$ only.
The result (\ref{preceding}) was later improved in several papers by sharpening the moment assumptions and relaxing the restrictions between $n$ and $p$. In terms of the relationship between $n$ and $p$, these results can be classified into the following categories:
\begin{itemize}
\item[(a).]  {\it Linear rate: $p \sim c n$ with $c$ being a constant}.
Li and Rosalsky (2006), Zhou (2007), Li, Liu and Rosalsky (2009) and Li, Qi and Rosalsky (2010) improved the moment conditions to make (\ref{preceding}) valid under  the condition $p/n \to c\in (0, 1).$

\item[(b).]  {\it Polynomial rate: $p = O(n^{\alpha})$ with $\alpha>0$ being a constant}. Liu, Lin and Shao (2008) showed that (\ref{preceding}) holds as $p\to \infty$ and $p = O(n^{\alpha})$ where $\alpha$ is a constant. That is, (\ref{preceding}) still holds when $n$ and $p$ are in the polynomial rates.

\item[(c).] {\it Sub-exponential rate: $\log p=o(n^{\alpha})$ with $0<\alpha\le 1/3$ being a constant}.  Motivated by applications in testing high-dimensional covariance structure and construction of compressed sensing, Cai and Jiang (2010) further extended the range of $p$ by considering the sub-exponential rate. It was shown that (\ref{preceding}) is also valid if $\log p=o(n^{\alpha})$ with $\alpha \in (0, 1/3]$ and the distribution of $x_{11}$ is well-behaved.  In particular, (\ref{preceding}) holds with $\alpha =1/3$ when $x_{ij}$'s are i.i.d. $N(0,1)$ random variables.
\end{itemize}

An interesting question is whether the limiting distribution (\ref{preceding}) holds for even higher dimensional case when $\log p$ is of order $n^{\alpha}$ with $\alpha >1/3$. This is a case of significant interest in high-dimensional data analysis and signal processing. For example, in the context of high-dimensional regression and classification, simulation studies about the distribution of $L_n$ were made in Cai and Lv (2007) and Fan and Lv (2008 and 2010).
In this paper we shall study the limiting laws of the coherence $L_n$ for a full range of the values of $p$. To make our technical analysis tractable, we focus on the setting where the columns $\bd{X}_i$ of the random matrix $\bd{X}$ follow a spherical distribution, which contains the normal distribution $N(\bd{0}, \sigma^2\bd{I}_n)$ as a special case. Motivated by the applications in statistics and signal processing mentioned earlier, we are especially interested in the ultra high dimensional case.   More specifically, we consider three different regimes:
\begin{itemize}
\item[(i).]  {\it the sub-exponential  case:} $\frac{1}{n}\log p \to 0$;

\item[(ii).] {\it the exponential case:} $\frac{1}{n}\log p \to \beta \in (0, \infty)$;

\item[(iii).] {\it the super-exponential case:} $\frac{1}{n}\log p \to\infty.$
\end{itemize}
Our results show that the limiting behavior of $L_n$ differs significantly in different regimes and exhibits interesting phase transition phenomena as the dimension $p$ grows as a function of $n$. To answer the question posed earlier, it is shown that $nL_n^2-4\log p + \log\log p$ converges to the limiting distribution given in (\ref{preceding})  if and only if  $\log p = o(n^{1/2})$. The phase transition in the limiting distribution first occurs with the case when $\log p$ is of order $n^{1/2}$.  In this transitional case, additional shift in the limiting distribution occurs. When the dimension $p$ further grows as a function of $n$, another transition occurs in the range when $\log p$ is of the same order as  $n$. In the sub-exponential case, $L_n$ converges to $0$ in probability. When $\log p \sim \beta n$ for some positive constant $\beta$, $L_n$ converges in probability to a constant strictly between 0 and 1, and the limiting distribution of $T_n = \log (1-L_n^2)$ is significantly different from that in the sub-exponential case. If $p$ is further increased to the super-exponential case, $L_n$ converges to 1 in probability and the limiting distribution of $T_n$ becomes the extreme value distribution without a shift.

There are also interesting differences between the limiting behaviors of $L_n$  and $\tilde L_n$. As shown in Cai and Jiang (2010), the limiting laws of $L_n$  and $\tilde L_n$ coincide with each other when $x_{ij}$'s are iid $N(0, 1)$ variables and $\log p = o(n^{1/3})$. Our results show that this remains true in the current setting for the sub-exponential and exponential cases, but not true for the super-exponential case.  It is interesting to contrast the results obtained in this paper with the results on $L_n$ and $\tilde L_n$ in the previous literature. The only known limiting distribution of $L_n$ and $\tilde L_n$ is given in (\ref{preceding}) and the best known result  in terms of the range of $p$ is $\log p = o(n^{1/3})$. In comparison,  our study significantly extends the knowledge on the limiting distributions of the coherence and shows the ``colorful" phase transition phenomena as the dimension $p$ increases.

The limiting laws of the coherence have immediate  applications in statistics and signal processing.
Testing the covariance structure of a high dimensional random variable is an important problem in statistical inference. A particularly interesting problem is to test for independence in the Gaussian case because many statistical procedures are built upon the assumptions of independence and normality of the observations. The limiting laws of the coherence derived in this paper can be used directly to construct a test for independence in the ultra high dimensional setting. In addition, the limiting laws can also be used for the construction of compressed sensing matrices. We shall discuss these applications in Section \ref{app.sec}.

Many sophisticated probabilistic tools have been used in the previous literature to study the limiting laws of the coherence.
For example, the Chen-Stein method, large deviation inequalities, and strong approximations were used to derive the results mentioned earlier in (a), (b) and (c). Yet there appears to be limitations to these methods. It is unclear (to us) whether these techniques can be easily adopted to derive the limiting distribution of $L_n$ when $\log p$ is of order $n^{\alpha}$ for $\alpha >1/3$ and answer the question posed earlier. See Remark \ref{honey} in Section \ref{soft} for further discussions. In this paper  a different technique is developed. Under the assumption that $\bd{X}_i$ in (\ref{mount}) has a spherical distribution, we first show a somewhat surprising result that the sample correlation coefficients $\{\rho_{ij}; \; 1\le i < j\le p\}$ are pairwise independent. We then apply the Chen-Stein method to the coherence $L_n=\max_{1\leq i <  j \leq p}|\rho_{ij}|$ by using the exact distribution of $\rho_{ij}$ and the pairwise-independence structure of $\rho_{ij}$. In addition, the exact distribution of $\rho_{ij}$ also leads to some interesting properties of $\rho_{ij}$ in the small sample cases: $\rho_{ij}$ has the {\it symmetric Bernoulli distribution} for $n=2$, that is, $P(\rho_{ij}=\pm 1)=1/2$; $\rho_{ij}^2$ follows the {\it Arcsine law} on $[0, 1]$ for $n=3;$  $\rho_{ij}$ follows the {\it uniform distribution} on $[-1, 1]$ for $n=4;$ and  $\rho_{ij}$  follows the {\it semi-circle law} for $n=5.$

The rest of the paper is organized as follows. Section \ref{limit.sec} studies the limiting laws of the coherence $L_n$  and $\tilde L_n$ of a random matrix in the high-dimensional setting under the three regimes. The interesting phase transition phenomena are  discussed in detail.
Section \ref{app.sec} considers two direct applications of the limiting laws derived in this paper to statistics and signal processing in the ultra high dimensional setting.  Section \ref{soft} discusses some of the interesting aspects of the techniques used in the derivations.
Connections and differences with other related work, for example, the relationship between the sample correlation coefficients and the angles between random vectors, are discussed in Section \ref{discussion.sec}. The main results are proved in Section \ref{proof.sec}.

\section{Limiting Laws of the Coherence}
\label{limit.sec}

In this section we study separately the limiting behaviors of the coherence $L_n$  and $\tilde L_n$ of an $n\times p$ random matrix $\bd{X}$ under the three regimes: $\frac{1}{n}\log p \to 0$, $\frac{1}{n}\log p \to \beta\in (0, \infty)$, and $\frac{1}{n}\log p \to\infty$. As mentioned before, we shall focus on the setting where the columns $\bd{X}_i$ of the random matrix $\bd{X}$ follow a spherical distribution.

\subsection{Limiting Laws of the Coherence $L_n$}
\lbl{suite}

A random vector $\bd{Y}\in \mathbb{R}^n$ is said to follow a {\it spherical distribution} if $\bd{O}\bd{Y}$ and $\bd{Y}$ have the same probability distribution for all $n\times n$ orthogonal matrix $\bd{O}.$ Examples of spherical distributions include:

\begin{itemize}
\item  the multivariate normal distribution $N(\bd{0}, \sigma^2\bd{I}_n)$ with $\sigma>0$;

\item the  normal scale-mixutre distribution $\sum_{k=1}^K \epsilon_k N(\bd{0}, \sigma_k^2\bd{I}_n)$ with the density function
\bea\lbl{multiCon}
\sum_{k=1}^K \epsilon_k \frac{1}{(2\pi\sigma_k^2)^{n/2}}\cdot \exp\left(-\frac{1}{2\sigma_k^2}\bd{y}^T\bd{y}\right)
\eea
where $\sigma_k > 0$, $\epsilon_k>0$, and $\sum_{k=1}^K \epsilon_k=1$;

\item  the multivariate $t$ distribution with $m$ degrees of freedom and density function
\bea\lbl{multiT}
\frac{\Gamma(\frac{m+n}{2})}{\Gamma(\frac{m}{2})(m\pi)^{n/2}}\cdot \Big(1+\frac{1}{m}\bd{y}^T\bd{y}\Big)^{(m+n)/2}
\eea
for $m\geq 1.$ The case $m=1$ corresponds to the multivariate Cauchy distribution.
\end{itemize}
See Muirhead (1982) for further discussions on spherical distributions.

Let $\bd{X}=(\bd{X}_1, \cdots, \bd{X}_p)=(x_{ij})_{n\times p}$ be an $n\times p$ random matrix. Throughout the rest of this paper, we shall assume:

\medskip \noindent
 {\bf Assumption (A)}:  the columns $\bd{X}_1, \cdots, \bd{X}_p$ are independent $n$-dimensional random vectors with a common spherical distribution (which may depend on $n$) and $P(\bd{X}_1=\bd{0})=0$.

\medskip
The condition $P(\bd{X}_1=\bd{0})=0$  is to ensure that the correlation coefficients are well defined.
Let $\rho_{ij}$ be the Pearson correlation coefficient of $\bd{X}_i$ and $\bd{X}_j$ for $1\leq i < j\leq p.$ Then, $\Psi_n:=(\rho_{ij})_{p\times p}$ is the correlation matrix of $\bd{X}$, and $L_n$ defined in (\ref{mount}), is the largest magnitude of the off-diagonal entries of the sample correlation matrix $\Psi_n.$

To make the statements of the limiting distributions uniform across different regimes, we shall state all the results in the main theorems in terms of $T_n=\log (1-L_n^2)$. We begin with the sub-exponential case.

\begin{theorem}[Sub-Exponential Case]
\lbl{kind}
Suppose $p=p_n$ satisfies $(\log p)/n \to 0$ as $n\to\infty$,  then under Assumption (A),

\bei
\item[\rm (i).] $L_n\to 0$ in probability as $n\to\infty.$

\item[\rm (ii).] Let $T_n=\log (1-L_n^2).$ Then, as $n\to\infty,$
\bea\lbl{river}
nT_n + 4\log p -\log \log p
\eea
 converges weakly to an extreme distribution with the distribution function $F(y)=1- e^{-Ke^{y/2}},\  y\in\mathbb{R}$ and $K=1/\sqrt{8\pi}.$
\eei
\end{theorem}

The following law of large numbers is a direct consequence of Theorem \ref{kind}.

\begin{coro}\lbl{sister} Assume the same conditions as in Theorem \ref{kind}, we have
\bea\lbl{wiggle}
\sqrt{\frac{n}{\log p}}L_n \to 2
\eea
in probability as $n\to\infty.$
\end{coro}
This result actually provides the convergence speed of $L_n \to 0$ stated in Theorem \ref{kind}(i). It is stronger than Theorem 2 in Cai and Jiang (2010), which shows  (\ref{wiggle}) holds if  $\log p =o(n^{1/3})$ and $x_{ij}$'s are i.i.d. $N(0,1)$  random variables.

Theorem \ref{kind} also shows an interesting phase transition phenomenon of the limiting behavior of the coherence $L_n$.

\begin{coro}[Transitional Case]
\lbl{brother}
Suppose $p=p_n$ satisfies $\lim_{n\to\infty}(\log p)/\sqrt{n}=\alpha\in [0, \infty)$, then under Assumption (A),
\bea\lbl{glass}
nL_n^2 - 4\log p +\log \log p
\eea
converges weakly to the distribution function $\exp\{-{1\over \sqrt{8\pi}} e^{-(y+8\alpha^2)/2}\}$,  $y\in \mathbb{R}$.
\end{coro}
As mentioned in the introduction, Cai and Jiang (2010) shows that $nL_n^2-4\log p + \log\log p$ converges weakly to an extreme distribution with distribution function given in (\ref{preceding}) when $\log p = o(n^{1/3})$ and  $x_{ij}$ are independent standard normal variables.  This is the best known result in the literature in terms of the range of $p$. Corollary \ref{brother} shows that (\ref{preceding}) holds if and only if  $\log p = o(n^{1/2})$ when $\bd{X}_1$ has a spherical distribution which includes the normal distribution $N(0, I_n)$ as a special case. This answers the question asked earlier in this paper.  Corollary \ref{brother} also shows that the limiting distribution of $L_n$ has a transitional phase between $(\log p)/\sqrt{n}\to 0$ and $(\log p)/\sqrt{n}\to \infty.$ In the transitional case when $(\log p)/\sqrt{n}\to \alpha\in (0, \infty)$, the limiting distribution of $nL_n^2 - 4\log p +\log \log p$ is shifted to the left  by $8\alpha^2.$


We now consider the exponential case.
\begin{theorem}[Exponential Case]
\lbl{jean}
Suppose $p=p_n$ satisfies $(\log p)/n \to \beta \in (0, \infty)$ as $n\to\infty$, then under Assumption (A),
\bei
\item[\rm (i).] $L_n\to \sqrt{1-e^{-4\beta}}$ in probability as $n\to\infty.$

\item[\rm (ii).] Let $T_n=\log (1-L_n^2).$ Then, as $n\to\infty,$
\bea\lbl{divide}
nT_n + 4\log p -\log \log p
\eea
 converges weakly to the distribution function
\bea\lbl{east}
F(y)=1- \exp\left\{-K(\beta) e^{(y+8\beta)/2}\right\},\  y\in\mathbb{R},\ \mbox{where}\    K(\beta)=\Big(\frac{\beta}{2\pi(1- e^{-4\beta})}\Big)^{1/2}.
\eea
\eei
\end{theorem}
Theorem \ref{jean} reveals the behavior of $L_n$ in the transitional case $(\log p)/n \to \beta$. In this case, the coherence $L_n$ converges in probability to a constant strictly between 0 and 1. Dividing  (\ref{divide}) by $n$, it is easy to see that
\[
\mbox{\rm $T_n\to -4\beta$ in probability as $n\to\infty$}
\]
since  $\lim_{n\to\infty}(\log p)/n = \beta \in (0, \infty).$ This is  also a direct consequence of  Theorem \ref{jean}(i).
Furthermore, it is trivially true that $1-e^{-4\beta}\sim 4 \beta$ as $\beta \to 0^+.$ Thus,
\beaa
\lim_{\beta \to 0^+}K(\beta)=\frac{1}{\sqrt{8\pi}},
\eeaa
which is exactly the value of $K$ in Theorem \ref{kind}. Thus, the limiting distribution $F(y)$ in Theorem \ref{jean}   as $\beta\to 0^+$ becomes the limiting distribution $F(y)$ in Theorem \ref{kind}. Heuristically, the sub-exponential case covered in Theorem \ref{kind} corresponds to the case ``$\beta=0$" in Theorem \ref{jean}. On the other hand, the exponential case of $(\log p)/n \to \beta\in (0, \infty)$ can also be viewed as a transitional phase between the sub-exponential and super-exponential cases.

Finally we turn to the super-exponential case where $(\log p)/n \to \infty$.

\begin{theorem}[Super-Exponential Case]
\lbl{learn}
Suppose $p=p_n$ satisfies $(\log p)/n \to \infty$ as $n\to\infty$. Let $T_n=\log (1-L_n^2).$ Then under Assumption (A),

\bei
\item[\rm (i).] $L_n \to 1$ in probability as $n\to\infty.$ Further, $\frac{n}{\log p} T_n  \to -4$ in probability as $n\to\infty.$

\item[\rm (ii).] As $n\to\infty,$
\bea\lbl{eve}
nT_n + \frac{4n}{n-2}\log p-\log n
\eea
 converges weakly to the distribution function $F(y)=1- e^{-Ke^{y/2}},\  y\in\mathbb{R}$ with $K=1/\sqrt{2\pi}.$
 \eei
\end{theorem}

The correction term of $nT_n$ in (\ref{eve}) is $\frac{4n}{n-2}\log p -\log n,$ which is different from the term ``$4\log p -\log \log p$" appeared in (\ref{river}) and (\ref{divide}). A reason  is that $T_n$ converges to a finite constant in probability  in Theorems \ref{kind} and \ref{jean}, whereas $T_n$ goes to $-\infty$  in probability  in Theorem \ref{learn}.
On the other hand, suppose $(\log p)/n \to \beta \in (0, \infty)$ and $\beta$ is large, then $\log n= \log\log p -\log \beta + o(1)$ and
\beaa
\frac{4n}{n-2}\log p = 4\log p + \frac{8}{n-2}\log p =4\log p + 8\beta + o(1)
\eeaa
as $n\to\infty.$ Consequently, the quantity in (\ref{eve}) becomes
\beaa
\left(nT_n + 4\log p -\log\log  p\right) +\ \mbox{constant} + o(1)
\eeaa
as $n\to\infty.$ The part in the parenthesis is the same as (\ref{river}) in Theorem \ref{kind} and (\ref{divide}) in Theorem \ref{jean}. This says that, heuristically, the results in Theorems \ref{kind}, \ref{jean} and \ref{learn} are consistent.

The formulation in the above theorems is in terms of  $T_n=\log (1-L_n^2)$ for uniformity. However,  one can easily change the expressions in terms of the coherence $L_n.$ For instance,
\beaa
P(n\log (1-L_n^2) + 4\log p -\log \log p\leq y)=P(L_n\geq  \sqrt{s_n}\,)
\eeaa
where
\bea\lbl{bang}
s_n:= 1 - \exp\left\{\frac{1}{n}(-4\log p +\log \log p +y)\right\}.
\eea

\subsection{Limiting Laws of $\tilde{L}_n$}
\lbl{silent}

We now study the limiting laws of the coherence $\tilde L_n$ defined in (\ref{corr'})  and  (\ref{mount3}). Note that under Assumption (A),  the columns $\bd{X}_1, \cdots, \bd{X}_p$ are independent $n$-dimensional random vectors with a common spherical distribution. By symmetry, it is easy to see that the mean $\mu = E \bd{X}_i = \bd{0}$ if it exists and hence
\be\lbl{corr'6}
\tilde{\rho}_{ij}=\frac{\bd{X}_i^T\bd{X}_j}{\|\bd{X}_i\|\cdot \|\bd{X}_j\|}\ \ \mbox{and}\ \ \tilde{L}_n=\max_{1\leq i <  j \leq p}|\tilde{\rho}_{ij}|.
\ee
As mentioned in the introduction,  Cai and Jiang (2010) showed that the limiting laws of $L_n$  and $\tilde L_n$ coincide with each other when $x_{ij}$'s are iid $N(0, 1)$ random variables and $\log p = o(n^{1/3})$. We shall show that this is still true in our current setting for the sub-exponential and exponential cases, but not true for the super-exponential case.
\begin{theorem}[Sub-Exponential \& Exponential Cases]
\lbl{kind4}
Under the same conditions, Theorems \ref{kind} and \ref{jean} and Corollaries \ref{sister} and \ref{brother} hold with $L_n$ replaced by $\tilde L_n$.
\end{theorem}

In  the super-exponential case,  the limiting behaviors  of  $\tilde L_n$ and  $L_n$ are different.

\begin{theorem}[Super-Exponential Case]
\lbl{learn4}
Suppose $p=p_n$ satisfies $(\log p)/n \to \infty$ as $n\to\infty$. Let $\tilde{T}_n=\log (1-\tilde{L}_n^2).$ Then under Assumption (A),
\bei
\item[\rm (i).] $\tilde{L}_n \to 1$ in probability as $n\to\infty.$ Further, $\frac{n}{\log p} \tilde{T}_n \to -4$ in probability as $n\to\infty.$

\item[\rm (ii).]  As $n\to\infty,$
\bea\lbl{eve4}
n\tilde{T}_n + \frac{4n}{n-1}\log p-\log n
\eea
 converges weakly to the distribution function $F(y)=1- e^{-Ke^{y/2}},\  y\in\mathbb{R}$ with $K=1/\sqrt{2\pi}.$
 \eei
\end{theorem}
Note the difference between (\ref{eve}) and (\ref{eve4}). When $(\log p)/n \to \infty$, the difference between $\frac{4n}{n-2}\log p$ and $\frac{4n}{n-1}\log p$ is not negligible.

\section{Applications}
\label{app.sec}

As mentioned in the introduction, the limiting laws of the coherence have a wide range of applications. Here we discuss briefly two immediate applications, one in high-dimensional statistics and another in signal processing. These applications were also discussed in Cai and Jiang (2010), but restricted to the Gaussian case with $\log p = o(n^{1/3})$. Here we extend to the more general spherical distributions and higher dimensions.

Testing the covariance structure of a distribution is an important problem in high dimensional statistical inference.
Let $\mathbf{Y}_{1},\ldots ,\mathbf{Y}_{n}$ be a random sample from a $p$-variate spherical distribution with covariance matrix $\Sigma _{p\times p}=(\sigma_{ij})$. We wish to test the hypotheses that $\Sigma$ is diagonal, i.e.,
\begin{equation}
\label{hypothesis}
H_0:  \sigma_{i,j} = 0 \;\mbox{for all $|i-j|\ge 1$}
\;\;\mbox{vs.} \;
H_a: \sigma_{i,j} \neq 0 \;\;\mbox{for some $|i-j|\ge 1$}.
\end{equation}
In the Gaussian case, this is the same as testing for independence.
The asymptotic distribution of  $L_{n}$ can be used to construct a convenient test
statistic for testing the hypotheses in (\ref{hypothesis}).
For example, in the case $\log p = o(n^{1/2})$, an approximate level $\alpha$ test is to
reject the null hypothesis $H_0$ whenever
\[
L_{n}^2 \ge n^{-1} \Big(4\log p - \log\log p - \log (8\pi) - 2 \log\log(1-\alpha)^{-1}\Big).
\]
It follows directly from Theorem \ref{kind} that the size of this test goes to $\alpha$ asymptotically as $n\to \infty$. This test was introduced in Cai and Jiang (2010) in the Gaussian case with the restriction that $\log p = o(n^{1/3})$.

Similarly, in the exponential (and sub-exponential) case,  set
\[
D_{n,p} =n T_n + 4\log p - \log\log p.
\]
Then Theorem \ref{jean} states that
\bea\lbl{sparrow}
P\left(D_{n,p} \le y\right) \goto 1-\exp\left(-K(\beta) e^{(y+8\beta)/2}\right),
\eea
where $K(\beta)=\Big(\frac{\beta}{2\pi(1- e^{-4\beta})}\Big)^{1/2}$. An approximate level $\alpha$ test for testing the hypotheses in (\ref{hypothesis}) can be obtained by rejecting the null hypothesis $H_0$ whenever
\[
D_{n,p} \le  2 \log\log(1-\alpha)^{-1} - 2 \log K(\beta) - 8 \beta.
\]
A test for the super-exponential case can also be constructed analogously by using the limiting distribution given in Theorem \ref{learn}.


Compressed sensing is an active and fast growing field in signal processing. See, e.g.,
Donoho (2006), Candes and Tao (2007), Bickel,
Ritov and Tsybakov (2009), Candes and Plan (2009), and Cai, Wang and Xu (2010a, b).  An important problem in compressed sensing is the construction of measurement matrices $\bd{X}_{n\times p}$
which enables the precise recovery of a sparse signal $\bd{\beta}$ from
linear measurements $\bd{y} = \bd{X}\bd{\beta}$  using an efficient recovery
algorithm. Such a measurement matrix $\bd{X}$ is typically randomly generated because it is difficult to construct deterministically.
The best known example is perhaps the $n\times p$ random matrix
$\bd{X}$ whose entries $x_{i,j}$ are iid normal variables
\be
\label{normal.CS}
x_{i,j}\stackrel{iid}{\sim} N(0, n^{-1}).
\ee
A commonly used condition is the mutual incoherence property (MIP) which requires the pairwise correlations among the column vectors of $\bd{X}$ to be small. Write $\bd{X}=(\bd{X}_1, \cdots, \bd{X}_p)=(x_{ij})_{n\times p}$ with $x_{ij}$ satisfying (\ref{normal.CS}) and let the coherence
$\tilde{L}_n=\max_{1\leq i <  j \leq p}|\tilde{\rho}_{ij}|$ be defined as in  (\ref{corr'})  and  (\ref{mount3}).
It has been shown that the condition
\be
\label{sharp.condition1}
(2k-1)\tilde{L}_n < 1
\ee
ensures the exact recovery of $k$-sparse signal $\beta$ in the noiseless case where $y=X\beta$ (see Donoho and Huo (2001) and Fuchs (2004)), and stable recovery of sparse signal in the
noisy case where
\[
y=X\beta + z.
\]
Here $z$ is an error vector, not necessarily random. See Cai, Wang and Xu (2010b).

The limiting laws derived in this paper can be
used to show how likely a random matrix satisfies the MIP condition
(\ref{sharp.condition1}). Take the sub-exponential case as an example.
By Theorem \ref{kind4}, as long as $(\log p)/n \to 0,$
\[
\tilde{L}_n \sim 2\sqrt{\log p \over n}.
\]
So in order for the MIP condition (\ref{sharp.condition1}) to hold,
roughly the sparsity $k$ should satisfy
\[
k < {1\over 4} \sqrt{n\over \log p}.
\]

\section{Technical Tool: Distribution of Correlation Coefficients}
\lbl{soft}

In this section we shall discuss the methodology used in our technical arguments. Sophisticated approximation methods such as the Chen-Stein method, large deviation bounds and strong approximations are the main ingredients in the proofs of the previous results in the literature including those given in Jiang (2004), Li and Rosalsky (2006), Zhou (2007),  Liu, Lin and Shao (2008), Li, Liu and Rosalsky (2009), Li, Qi and Rosalsky (2010), and Cai and Jiang (2010).
Though these technical tools work well for the cases when the dimension $p$ is not ultra high, it is far from clear to us whether/how these same tools can be used to derive the limiting distributions of the coherence $L_n$ for the three regimes considered in Section \ref{limit.sec}.

In this paper, a different approach is developed to derive the limiting distributions of $L_n$.  Assuming the $\bd{X}_i$'s have the spherical distribution, we find an interesting and useful property of the  correlation coefficients $\{\rho_{ij};\, 1\leq i < j\leq p\}$ and $\{\tilde \rho_{ij};\, 1\leq i < j\leq p\}$ given below.

\begin{lemma}
\lbl{and}
Let $n\geq 3$. Under Assumption (A), the Pearson correlation coefficient $\{\rho_{ij};\, 1\leq i < j\leq p\}$ are pairwise independent and identically distributed with  density function
\bea\lbl{meeting}
f(\rho) =\frac{1}{\sqrt{\pi}}\frac{\Gamma(\frac{n-1}{2})}{\Gamma(\frac{n-2}{2})}\cdot(1-\rho^2)^{\frac{n-4}{2}},\ \ \ |\rho|<1.
\eea
Similarly,  $\{\tilde{\rho}_{ij};\, 1\leq i < j\leq p\}$ are pairwise independent and identically distributed with  density
\bea\lbl{meeting'}
g(\rho) =\frac{1}{\sqrt{\pi}}\frac{\Gamma(\frac{n}{2})}{\Gamma(\frac{n-1}{2})}\cdot(1-\rho^2)^{\frac{n-3}{2}},\ \ \ |\rho|<1.
\eea
\end{lemma}

Note that the only difference between (\ref{meeting}) and (\ref{meeting'}) is the ``degree of freedom": replacing $n$ in (\ref{meeting'}) with $n-1$, one gets (\ref{meeting}). This is not difficult to understand by noting the definition of $\rho_{ij}=\frac{(\bd{X}_i-\overline{\bd{X}}_i)^T(\bd{X}_j-\overline{\bd{X}}_j)}{\|\bd{X}_i-\overline{\bd{X}}_i\|\cdot \|\bd{X}_j-\overline{\bd{X}}_j\|}.$ Heuristically, by subtracting $\overline{\bd{X}}_i$ from $\bd{X}_i$, the distribution of $\rho_{ij}$ becomes one degree less than that of $\tilde{\rho}_{ij}=\frac{\bd{X}_i^T\bd{X}_j}{\|\bd{X}_i\|\cdot \|\bd{X}_j\|}$.

Although $\{\rho_{ij};\, 1\leq i < j\leq p\}$ are pairwise independent, they are not mutually independent. In fact, recalling $\Psi=\Psi_n=(\rho_{ij})_{p\times p},$ the probability density function of $\Psi$ is given by
\bea\lbl{minus}
h(\Psi)=B_{n,p}\cdot (\det(\Psi))^{(n-p-2)/2}\ \ \ \  \ \ \ (|\rho_{ij}|<1,\ i<j)
\eea
for $1\le p < n$, where $B_{n,p}$ is an (explicit) normalizing constant, see p.148 from Muirhead (1982). Obviously, $h(\Psi)$ is not a product of functions of individual $\rho_{ij}$'s, the entries of  $\Psi$, hence $\{\rho_{ij};\, 1\leq i < j\leq p\}$ are not independent.

Lemma \ref{and} also yields the following interesting results on the distribution of the correlation coefficients $\rho_{ij}$ in the small sample cases. The verification is given in Section \ref{proof.sec}.

\begin{coro}
\lbl{touch}
Under Assumption (A),  the following holds for all $1\leq i < j\leq p$.
\bei
\item[\rm (i).] When $n=2,$ $\rho_{ij}$ has the symmetric Bernoulli distribution, i.e., $P(\rho_{ij}=\pm 1)=1/2.$

\item[\rm (ii).] When $n=3,$ $\rho_{ij}$ has the density $f(\rho)=\frac{1}{\pi}\frac{1}{\sqrt{1-\rho^2}}$ on $(-1, 1).$
That is, $\rho_{ij}^2$ follows the arcsine law on $[0, 1]$.

\item[\rm (iii).] When $n=4,$  $\rho_{ij}$
follows the uniform distribution on $[-1, 1].$

\item[\rm (iv).] When $n=5,$ $\rho_{ij}$ has the density $f(\rho)=\frac{2}{\pi}\sqrt{1-\rho^2}$ for $|\rho|\leq 1.$ That is, $\rho_{ij}$ follows the semi-circle law.
\eei
\end{coro}
Lemma \ref{and} provides a major technical tool for the proof of the main results. The starting step in the proofs of our theorems is the Chen-Stein method (Lemma \ref{stein}) which requires the evaluation of  two quantities: $P(\rho_{ij}\geq C)$ and $P(\rho_{ij}\geq C,\rho_{kl}\geq C).$ By using the explicit density expression in (\ref{meeting}), we are able to evaluate the first probability precisely.  The pairwise independence stated in Lemma \ref{and} yields $P(\rho_{ij}\geq C,\rho_{kl}\geq C)=P(\rho_{ij}\geq C)^2$ for $\{i, j\} \ne \{k, l\}.$ In other words, the evaluation of the second quantity is reduced to the study of the first one. This greatly simplifies some of the technical arguments.


\begin{remark}\lbl{honey}
{\rm
Equation (\ref{meeting}) yields directly that $W_n:=\sqrt{n}\rho_{12}$ has the density function 
\beaa
f_n(w)=\frac{1}{\sqrt{n}}\cdot\frac{1}{\sqrt{\pi}}\frac{\Gamma(\frac{n-1}{2})}{\Gamma(\frac{n-2}{2})}\cdot
\left(1-\frac{w^2}{n}\right)^{\frac{n-4}{2}}\to \frac{1}{\sqrt{2\pi}}e^{-w^2/2}
\eeaa
as $n\to\infty$ for all $w\in \mathbb{R},$ where the fact  that $\Gamma(\frac{n-1}{2})/\Gamma(\frac{n-2}{2}) \sim \sqrt{n/2}$ as $n\to\infty$ (see (\ref{white})) is used. This shows that $W_n$ converges to $N(0,1)$ in distribution as $n\to\infty.$ Set $(x_{ij})_{n\times p}:=(\bd{X}_1, \cdots, \bd{X}_p).$ Assuming that $x_{ij}$'s are i.i.d. with an unknown distribution but with suitable moment conditions, say, $|x_{12}|$ is bounded, it can be shown easily that $\sqrt{n}\rho_{12}$ converges to $N(0,1)$ by using the standard central limit theorem for i.i.d. random variables and the Slusky theorem. However, the convergence speed is hard to be captured well enough so that $L_n$ in (\ref{mount}) is understood clearly when $p$ is much larger than $n.$ The best known result is that (\ref{preceding}) holds for  $\log p=o(n^{\alpha})$ with $\alpha= 1/3$ in  Cai and Jiang (2010). Here, with the understanding of the pairwise independence among $\{\rho_{ij};\, 1\leq i<j\leq p\}$ and the exact distribution of $\rho_{ij}$ we are able to get the limiting distribution of $L_n$ for the full range of  the values of $p$ and to fully characterize the phase transition phenomena in the limiting behaviors of the coherence (Theorems \ref{kind}, \ref{jean} and \ref{learn} and the corresponding corollaries).
}
\end{remark}

\section{Discussions}
\label{discussion.sec}

The present paper was inspired by the applications in high-dimensional statistics and signal processing in which the dimension $p$ is often desired to be as high as possible as a function of $n$. All the known results on the coherence $L_n$ are restricted to the cases where the dimension $p$ is either linear, polynomial or at most sub-exponential in $n$. In comparison,  we give in this paper a complete characterization of the limiting distribution of $L_n$ for the full range of $p$ including
the sub-exponential  case $\frac{1}{n}\log p \to 0$,  the exponential case  $\frac{1}{n}\log p \to \beta\in (0, \infty)$, and the super-exponential case $\frac{1}{n}\log p \to\infty.$  Our results show interesting phase transition phenomena in the limiting distributions of the coherence when the dimension $p$ grows as a function of $n$. Over the full range of values of $p$, phase transition of the limiting behavior of $L_n$ occurs twice: when $\log p$ is of order $n^{1/2}$ and when $\log p$ is of order $n$. These results also show that the standard limiting distribution (\ref{preceding}) known in the literature holds if and only if $\log p = o(n^{1/2})$ when the columns have a spherical distribution which includes the commonly considered i.i.d. normal setting as a special case.

Previous results on the coherence $L_n$ focus on the  case where the entries $x_{ij}$ of the random matrix $\bd{X}$ are i.i.d. under certain moment conditions. See the references mentioned in (a), (b) and (c) in the introduction.
In this paper, we assume the columns of $\bd{X}=(x_{ij})_{n\times p}$ to be i.i.d. with a spherical distribution. The spherical distribution assumption are more special than the non-specified distributions with certain moment conditions considered in the previous literature. On the other hand,  the entries of a vector with a spherical distribution do not have to be independent (see, e.g., the normal scale-mixture distribution in (\ref{multiCon}) and the multivariate $t$-distribution in (\ref{multiT})). In this sense, our work relaxes the independence assumption among the entries $x_{ij}$.
Under the assumption of spherical distributions, we are able to show that the sample correlation coefficients are pairwise independent and then use the exact distribution and the pairwise-independence structure of the sample correlation coefficients as a major technical tool in the derivation of the limiting distributions.

There are interesting connections between sample correlation coefficients and angles between random vectors.  Let $\bold{a}\in \mathbb{R}^n$ be a deterministic vector with $\|\bold{a}\|=1.$ Let $\bold{X}_1\in \mathbb{R}^n$ be a random vector with a spherical distribution satisfying $P(\bd{X}_1 =\bd{0})=0.$ Relating Theorem 1.5.7(i) and (5) on page 147 in Muirhead (1982), it can be seen that $W=\frac{\bold{a}^T\bold{X}_1}{\|\bold{X}_1\|}$ has the same distribution as the one given in (\ref{meeting'}). Note that $\frac{\bold{X}_1}{\|\bold{X}_1\|}$ has the uniform distribution over the unit sphere in $\mathbb{R}^n$, and hence $W$ is the cosine of the angle between a fixed unit vector $\bold{a}$ and a random vector with the uniform distribution on the unit sphere. Similar to Corollary \ref{touch}, the following holds.
\bei
\item[\rm(i).]  If $n=2,$ then the cosine of the angle has the probability density function $f(\rho)=\frac{1}{\pi}\frac{1}{\sqrt{1-\rho^2}}.$
That is, the square of the cosine follows the Arcsine law on $[0, 1]$.

\item[\rm (ii).] If $n=3,$ then the cosine of the angle
follows the uniform distribution on $[-1, 1].$

\item[\rm (iii).] If $n=4,$  then  the cosine of the angle has the probability density function $f(\rho)=\frac{2}{\pi}\sqrt{1-\rho^2}$ for $|\rho|\leq 1.$ That is, $\rho_{ij}$ follows the semi-circle law.
\eei

The semi-circle law is perhaps best known in random matrix theory as the limit of the empirical distribution of the eigenvalues of an $n\times n$ Wigner random matrix as $n\to \infty$. See, e.g., Wigner  (1958).  It seems not so common to see a random variable to satisfy the semi-circle law in practice. It is interesting to see the semi-circle law here as the exact distribution of the correlation coefficient and the cosine of the angle between two random vectors in  Corollary \ref{touch}(iv) and (iii) above.

\section{Proofs}
\lbl{proof.sec}
In this section we prove the main results of the paper. We shall write $p$ for $p_n$ if there is no confusion. We begin by proving Lemma \ref{and} on the distributions of  the correlation coefficients. We then collect and prove a few additional technical results before giving the proofs of the main theorems.




\subsection{Technical Results}\lbl{suite'}

The following lemma is needed for the proof of Lemma \ref{and}.

\begin{lemma}\label{measure0}
Let $\bd{X}$ be an $n$-dimensional random vector with a spherical distribution and $P(\bd{X}=\bd{0})=0$.  Let $\bd{1}=(1,\cdots, 1)^T\in \mathbb{R}^n$ and $\{1\}=\{k\bd{1};\, k\in \mathbb{R}\},$ the span of $\bd{1}.$ Then $P(\bd{X}\in \{1\})=0$.
\end{lemma}

\noindent\textbf{Proof.} Since $P(\bd{X}=\bd{0})=0$, we know $\bd{Y}:=\frac{\bd{X}}{\|\bd{X}\|}$ is well-defined. By definition, $\bd{O}\bd{X} \overset{P}{=} \bd{X}$ for any orthogonal matrix $\bd{O}$, then 
\beaa
\bd{O}\bd{Y}=\frac{\bd{O}\bd{X}}{\|\bd{O}\bd{X}\|} \overset{P}{=}\frac{\bd{X}}{\|\bd{X}\|}=\bd{Y}.
\eeaa
That is, the probability measure generated by $\bd{Y}$ is an orthogonal-invariant measure on the unit sphere $S^{n-1} \subset\mathbb{R}^n.$ Since the Haar probability measure, as the distribution on the unit sphere with the orthogonal-invariant property, is unique, it follows that  $\bd{Y}$ must have the uniform distribution on the unit sphere in $\mathbb{R}^n.$ In particular, $P(\bd{Y}=y)=0$ for any $y\in S^{n-1}.$ Let $A=\{\bd{X}\in \{1\}\backslash \{\bd{0}\}\}$ and $y_0={n}^{-1/2}(1,\cdots,1)^T\in S^{n-1}.$ Notice
$A\subset \{\bd{Y}=y_0\ \mbox{or}\, -y_0\}.$ It follows that $P(\bd{X}\in \{1\})=P(A)\leq P(\bd{Y}=y_0)+ P(\bd{Y}=-y_0) =0.$ \ \ \ \ \ \ \ $\blacksquare$\\

\noindent\textbf{Proof of Lemma \ref{and}}. Recall that $\bd{X}_1, \cdots, \bd{X}_p$ are independent and  $\rho_{ij}$ is the Pearson correlation coefficient of $\bd{X}_i$ and $\bd{X}_j$ for $1\leq i < j\leq p.$ Given $i<j$ and $k<l$ with $(i,j)\ne (k,l).$ It is easy to see that $\rho_{ij}$ and $\rho_{kl}$ are independent if $\{i,j\}\bigcap \{k,l\}=\varnothing.$ Thus, to finish the proof, it enough to prove the following:
\bea\lbl{voice}
& & \mbox{Let}\ \{\bd{U}, \bd{V}, \bd{W}\}\, \mbox{be i.i.d with an}\, n\mbox{-dimensional spherical distribution and}\ P(\bd{U}=\bd{0})=0.  \nonumber\\
& & \mbox{Then}\ \rho_{\bd{U},\bd{V}}\ \mbox{and}\ \rho_{\bd{U},\bd{W}}\ \mbox{are i.i.d. with the density function given in (\ref{meeting}).}
\eea
By Lemma \ref{measure0}, $P(\bd{U}\in \{1\})=P(\bd{V}\in \{1\})=P(\bd{W}\in \{1\})=0.$ Then,  $\rho_{\bd{U},\bd{V}}$ and $\rho_{\bd{U},\bd{W}}$ have the same  probability density function $f(\rho)$ by  (5) on p. 147 from Muirhead (1982). To show the independence, we need to prove
\bea\lbl{contest}
E[g(\rho_{\bd{U},\bd{V}})\cdot h(\rho_{\bd{U}, \bd{W}})]=Eg(\rho_{\bd{U},\bd{V}})\cdot Eh(\rho_{\bd{U},\bd{W}})
\eea
for any bounded and measurable functions $g(x)$ and $h(x).$ Since $\bd{U}$, $\bd{V}$ and $\bd{W}$ are independent,
\bea\lbl{look}
E[g(\rho_{\bd{U},\bd{V}})\cdot h(\rho_{\bd{U}, \bd{W}})]&=&E\Big\{E[g(\rho_{\bd{U},\bd{V}})\cdot h(\rho_{\bd{U},\bd{W}})|\bd{U}]\Big\}\nonumber\\
&=& E\Big\{E[g(\rho_{\bd{U},\bd{V}})|\bd{U}]\cdot E[h(\rho_{\bd{U},\bd{W}})|\bd{U}]\Big\}.
\eea
Write $\bd{V}=(V_1, \cdots, V_n)^T\in \mathbb{R}^n$ and $\bar{\bd{V}}=\frac{1}{n}\sum_{i=1}^nV_i$. For any numbers $u_1, \cdots, u_n$ such that at least two of them are not identical,  Theorem 5.1.1 and (5) on p. 147 from Muirhead (1982) say that
\beaa
\rho_{\bd{u}, \bd{V}}=\frac{\sum_{i=1}^n(u_i-\bar{\bd{u}})(V_i - \bar{\bd{V}})}
{\sqrt{\sum_{i=1}^n(u_i-\bar{\bd{u}})^2\cdot \sum_{i=1}^n(V_i-\bar{\bd{V}})^2}}
\eeaa
has the probability density function $f(\rho)$ as in (\ref{meeting}), where $\bd{u}=(u_1, \cdots, u_n)^T$ and  $\bar{\bd{u}}=\frac{1}{n}\sum_{i=1}^nu_i$ (see also Kariya and Eaton (1977) for this).
In other words, given $\bd{U}$, the probability distribution of $\rho_{\bd{U},\bd{V}}$ does not depend on the value of $\bd{U}.$  Let $\bd{U}=(U_1, \cdots, U_n)^T.$ Evidently,  $P(U_1=\cdots =U_n)=P(\bd{U} \in \{1\})=0$. Thus,
\[
 E[g(\rho_{\bd{U},\bd{V}})|\bd{U}]=\int_{|\rho|\leq 1}g(\rho)f(\rho)\,d\rho=E g(\rho_{\bd{U},\bd{V}})
 \]
 and
\[
E[h(\rho_{\bd{U},\bd{W}})|\bd{U}]=Eh(\rho_{\bd{U},\bd{W}})
\]
since $\rho_{\bd{U},\bd{V}}$ and $\rho_{\bd{U},\bd{W}}$ have the same probability density function $f(\rho)$ as in (\ref{meeting}). These and (\ref{look}) conclude (\ref{contest}).

We now turn to study $\tilde{\rho}_{ij}$. Given $1\leq i < j \leq p.$ Then $\bd{\alpha}:=\frac{\bd{X}_j}{\|\bd{X}_j\|}$ is a unit vector and is independent of $\bd{X}_i.$ Further, $\tilde{\rho}_{ij}=\frac{\bd{\alpha}^T\bd{X}_i}{\|\bd{X}_i\|}.$ It then follows from Theorem 1.5.7(i) and the argument for (5) on p.147 of Muirhead (1982) that
$\tilde{\rho}_{ij}$ has the probability density function $f(\rho)$ as in (\ref{meeting'}). The proof for the pairwise independence among $\{\tilde{\rho}_{ij};\, 1\leq i < j\leq p\}$ is the same as that for  the $\rho_{ij}$'s. \ \ \ \ \ \ \ \ \ $\blacksquare$\\

\noindent\textbf{Proof of Corollary \ref{touch}}. Taking $n=3,4,5,$ respectively, in Lemma \ref{and}, we easily have (ii), (iii) and (iv). Now we check (i).

Let $\bd{X}_1=(\xi_1, \eta_1)^T\in \mathbb{R}^2$ and $\bd{X}_2=(\xi_2, \eta_2)^T\in \mathbb{R}^2.$ It is easy to see
\bea\lbl{aspects}
\rho_{1\,2}=\frac{\xi_1-\eta_1}{|\xi_1-\eta_1|}\cdot \frac{\xi_2-\eta_2}{|\xi_2-\eta_2|}.
\eea
First, Assumption (A) and Lemma \ref{measure0} imply $P(\xi_i=\eta_i)=0$ for $i=1,2$. Since $\bd{X}_1$ has a spherical distribution, we know that $\bd{A}\bd{X}_1$ and $\bd{X}_1$ have the same distribution for any $\bd{A}=\mbox{diag}(\epsilon_1, \epsilon_2)$ with $\epsilon_i=\pm 1$, $i=1,2.$ This implies $\bd{X}_1$ is symmetric, and hence $\xi_1-\eta_1$ is symmetric. Consequently, $\frac{\xi_1-\eta_1}{|\xi_1-\eta_1|}$ takes value $\pm 1$ with probability $1/2$ each. The same is true for $\frac{\xi_2-\eta_2}{|\xi_2-\eta_2|}.$ By (\ref{aspects}) and the  independence between $\bd{X}_1$ and $\bd{X}_2,$ we conclude that $P(\rho_{1\,2}=\pm 1)=1/2.$\ \ \ \ \ \ \ \ $\blacksquare$

\begin{lemma}\lbl{also} Let $t=t_m\in(0,1)$ satisfy $mt_m^2\to \infty$ as $m\to\infty.$ Then
\beaa
\int_{t}^1(1-x^2)^{m/2}\,dx =\frac{1}{m t}\left(1-t^2\right)^{(m+2)/2}(1+o(1))
\eeaa
as $m\to\infty.$
\end{lemma}
\textbf{Proof}. Set $y=x^2$ for $x>0.$ Then $x=\sqrt{y}$ and
\bea
I_m:=\int_{t}^1(1-x^2)^{m/2}\,dx
&= &\frac{1}{2}\int_{t^2}^1\frac{1}{\sqrt{y}}(1-y)^{m/2}\,dy\lbl{buddy}\\
&= & -\frac{1}{m+2}\int_{t^2}^1\frac{1}{\sqrt{y}}\Big[(1-y)^{(m+2)/2}\Big]'\,dy.\nonumber
\eea
By integration by parts,
\bea
I_m &= & -\frac{1}{m+2}\frac{1}{\sqrt{y}}(1-y)^{(m+2)/2}\Big|_{t^2}^1
-\frac{1}{2(m+2)}\int_{t^2}^1\frac{1}{y^{3/2}}(1-y)^{(m+2)/2}\,dy\nonumber\\
& = & \frac{1}{(m+2)t}\left(1-t^2\right)^{(m+2)/2} -\frac{1}{m+2}\cdot\frac{1}{2}\int_{t^2}^1\frac{1}{\sqrt{y}}(1-y)^{m/2}\cdot \frac{1-y}{y}\,dy.\lbl{turk}
\eea
Note that $0\leq \frac{1-y}{y} \leq \frac{1}{t^2}$ for all $[t^2, 1].$ By the second equality in  (\ref{buddy}),
\beaa
0<\frac{1}{m+2}\cdot\frac{1}{2}\int_{t^2}^1\frac{1}{\sqrt{y}}(1-y)^{m/2}\cdot \frac{1-y}{y}\,dy \leq \frac{1}{mt^2} I_m.
\eeaa
This and (\ref{turk}) conclude that
\beaa
\frac{1}{(m+2)t}\left(1-t^2\right)^{(m+2)/2}  - \frac{1}{mt^2} I_m\leq I_m \leq \frac{1}{(m+2)t}\left(1-t^2\right)^{(m+2)/2}.
\eeaa
Solving the first inequality on $I_m$, we have
\beaa
\Big(1+\frac{1}{mt^2}\Big)^{-1}\frac{1}{(m+2)t}\left(1-t^2\right)^{(m+2)/2}   \leq I_m \leq \frac{1}{(m+2)t}\left(1-t^2\right)^{(m+2)/2}.
\eeaa
By the given condition that $mt^2=mt_m^2\to\infty,$ we arrive at
\beaa
 I_m = \frac{1}{(m+2)t}\left(1-t^2\right)^{(m+2)/2}(1+o(1)) = \frac{1}{m t}\left(1-t^2\right)^{(m+2)/2}(1+o(1))
\eeaa
as $m\to\infty.$\ \ \ \ \ \ \ \ \ \ \ $\blacksquare$\\

The following Poisson approximation result
is essentially a special case of
Theorem $1$ from Arratia et al. (1989).
\begin{lemma}\label{stein} Let $I$ be an index set and $\{B_{\alpha}, \alpha\in I\}$ be a set of subsets of $I,$ that is, $B_{\alpha}\subset I$ for each $\alpha \in I.$  Let also $\{\eta_{\alpha}, \alpha\in I\}$ be random variables. For a given $t\in \mathbb{R},$ set $\lambda=\sum_{\alpha\in I}P(\eta_{\alpha}>t).$ Then
\beaa
|P(\max_{\alpha \in I}\eta_{\alpha} \leq t)-e^{-\lambda}| \leq (1\wedge \lambda^{-1})(b_1+b_2+b_3)
\eeaa
where
\beaa
& & b_1=\sum_{\alpha \in I}\sum_{\beta \in B_{\alpha}}P(\eta_{\alpha} >t)P(\eta_{\beta} >t),\ \
 b_2=\sum_{\alpha \in I}\sum_{\alpha\ne \beta \in B_{\alpha}}P(\eta_{\alpha} >t, \eta_{\beta} >t),\\
 & & b_3=\sum_{\alpha \in I}E|P(\eta_{\alpha} >t|\sigma(\eta_{\beta}, \beta \notin B_{\alpha})) - P(\eta_{\alpha} >t)|,
\eeaa
and $\sigma(\eta_{\beta}, \beta \notin B_{\alpha})$ is the $\sigma$-algebra generated by $\{\eta_{\beta}, \beta \notin B_{\alpha}\}.$
In particular, if $\eta_{\alpha}$ is independent of $\{\eta_{\beta}, \beta \notin B_{\alpha}\}$ for each $\alpha,$ then $b_3=0.$
\end{lemma}

\begin{lemma}\lbl{childish} Let $L_n$ be as in (\ref{mount}) and Assumption (A) hold. For $\{t_n\in [0,1];\, n\geq 1\}$, set
\beaa
h_n = \frac{n^{1/2}p^2}{\sqrt{2\pi}}
\int_{t_n}^1(1-x^2)^{\frac{n-4}{2}}\,dx,\ n\geq 1.
\eeaa
If $\lim_{n\to\infty}h_n=\lambda\in [0, \infty),$ then $\lim_{n\to\infty}P(L_n\leq t_n) = e^{-\lambda}.$
\end{lemma}
\textbf{Proof}. For brevity of notation, we sometimes write $t=t_n$ if there is no confusion. First, take $I=\{(i,j);\ 1\leq i< j \leq p\}.$ For $u =(i,j) \in I,$ set $B_{u}=\{(k,l)\in I;\ \mbox{one of}\ k\ \mbox{and}\ l =i\ \mbox{or}\ j,\ \mbox{but}\ (k,l)\ne u\},\ \eta_{u}=|\rho_{ij}|$ and $A_{u}=A_{ij}=\{|\rho_{ij}| > t\}.$  By the i.i.d. assumption on $\bd{X}_1, \cdots, \bd{X}_p$ and Lemma \ref{stein},
\bea\lbl{season}
|P(L_n\leq t) - e^{-\lambda_n}| \leq b_{1,n}+b_{2,n}
\eea
 where
\bea\lbl{sheep}
\lambda_n=\frac{p(p-1)}{2}P(A_{12})
\eea
 and
\beaa
b_{1,n}\leq 2p^3P(A_{12})^2\ \mbox{and}\ b_{2,n} \leq 2p^3P(A_{12}A_{13}).
\eeaa
By Lemma \ref{and}, $A_{12}$ and $A_{13}$ are independent events with the same probability. Thus, from (\ref{sheep}),
\bea\lbl{English}
b_{1,n} \wedge b_{2,n} \leq 2p^3P(A_{12})^2 \leq \frac{8p\lambda_n^2}{(p-1)^2}\leq \frac{32\lambda_n^2}{p}
\eea
for all $p\geq 2.$ Now we compute $P(A_{12}).$ In fact, by Lemma \ref{and} again,
\beaa
P(A_{12})=\int_{1>|x|>t}f(x)\,dx &= & \frac{1}{\sqrt{\pi}}\frac{\Gamma(\frac{n-1}{2})}{\Gamma(\frac{n-2}{2})}
\int_{1>|x|>t}(1-x^2)^{\frac{n-4}{2}}\,dx\\
& = & \frac{2}{\sqrt{\pi}}\frac{\Gamma(\frac{n-1}{2})}{\Gamma(\frac{n-2}{2})}
\int_{t}^1(1-x^2)^{\frac{n-4}{2}}\,dx.
\eeaa
Recalling the Stirling formula (see, e.g., p.368 from  Gamelin (2001) or (37) on
p.204 from Ahlfors (1979)):
\begin{eqnarray*}
\log\Gamma(z)=z\log z - z -\frac{1}{2}\log z+ \log \sqrt{2\pi}
 +O\left(\frac{1}{x}\right)
\end{eqnarray*}
as $x=\mbox{Re}\,(z)\to \infty,$ it is easy to verify that
\bea\lbl{white}
\frac{\Gamma(\frac{n-1}{2})}{\Gamma(\frac{n-2}{2})} \sim \sqrt{\frac{n}{2}}
\eea
as $n\to\infty.$ Thus,
\beaa
P(A_{12}) \sim \frac{2n^{1/2}}{\sqrt{2\pi}}
\int_{t}^1(1-x^2)^{\frac{n-4}{2}}\,dx
\eeaa
as $n\to\infty.$ From (\ref{sheep}), we know
\beaa
\lambda_n \sim \frac{n^{1/2}p^2}{\sqrt{2\pi}}
\int_{t}^1(1-x^2)^{\frac{n-4}{2}}\,dx=h_n
\eeaa
as $n\to\infty.$  Finally, by (\ref{season}) and (\ref{English}), we know
\beaa
\lim_{n\to\infty}P(L_n\leq t) = e^{-\lambda}\ \ \mbox{if}\ \ \lim_{n\to\infty}h_n=\lambda \in [0, \infty).\ \ \ \ \ \ \blacksquare
\eeaa

\subsection{Proofs for Results on $L_n$ in Section \ref{suite}}\lbl{suite'''}




\noindent\textbf{Proof of Theorem \ref{kind}}. (i). Assume (ii) of the theorem holds. Since $(\log p)/n \to 0$ as $n\to\infty,$  dividing (\ref{river}) by $n$, we see that $\log (1- L_n^2) \to 0$ in probability, or equivalently, $L_n\to 0$ in probability as $n\to\infty.$

(ii). The proof here does not rely on the conclusion in (i). We claim that
\bea\lbl{watch}
(n-2)\log (1-L_n^2) + 4\log p -\log\log p
\eea
 converges weakly to the distribution function $F(y)=1- e^{-K e^{y/2}},\  y\in\mathbb{R}.$ Once this holds, using the condition that $\log p=o(n)$ and the same argument as in (i) above, we have $\log (1-L_n^2)\to 0$ in probability as $n\to\infty.$ Then by the Slusky lemma,
 \beaa
 n\log (1-L_n^2) + 4\log p -\log\log p
 \eeaa
converges weakly to the distribution function $F(y)=1- e^{-K e^{y/2}},\  y\in\mathbb{R}.$ We then obtain (\ref{river}). Now we prove (\ref{watch}).

Fix $y\in \mathbb{R}.$ Let $N=n-2$ and $t=t_n\in [0,1)$ such that
\bea\lbl{fish}
\log (1-t^2) =\frac{-4\log p +\log\log p +y }{N} \wedge 0.
\eea
From (\ref{fish}) and the assumption $\log p =o(n)$, we have that $t_n\to 0^+$ as $n\to\infty,$ and hence $\log (1-t^2) \sim -t^2$.  Thus, (\ref{fish}) implies
\bea\lbl{view}
t \sim \Big(\frac{4\log p - \log\log p -y}{N}\Big)^{1/2} \sim \frac{2\sqrt{\log p}}{\sqrt{N}}\ \mbox{and}\ Nt_n^2 \to\infty
\eea
as $n\to\infty.$ By (\ref{fish}) again,
\bea\lbl{China}
 P((n-2)\log (1-L_n^2) + 4\log p -\log \log p \geq y)=P(L_n\leq t)
\eea
as $n$ is large enough. Now let's compute $h_n$ in Lemma \ref{childish} for $\lim_{n\to\infty}P(L_n\leq t)$. Recall
\bea\lbl{user}
h_n = \frac{n^{1/2}p^2}{\sqrt{2\pi}}
\int_{t}^1(1-x^2)^{\frac{n-4}{2}}\,dx.
\eea
 From Lemma \ref{also} and the second assertion in (\ref{view}),
\beaa
n^{1/2}p^2
\int_{t}^1(1-x^2)^{(n-4)/2}\,dx &\sim & \frac{n^{1/2}p^2}{nt}\left(1-t^2\right)^{(n-2)/2}\\
& \sim & \frac{p^2}{\sqrt{N}\,t}\left(1-t^2\right)^{N/2}
\eeaa
as $n\to\infty.$ This joint with (\ref{fish}) and the first assertion in (\ref{view}) gives
\beaa
\frac{p^2}{\sqrt{N}\,t}\left(1-t^2\right)^{N/2} \sim \frac{p^2}{2\sqrt{\log p}}\cdot \exp\Big\{\frac{-4\log p + \log\log p +y}{N}\cdot \frac{N}{2}\Big\}=\frac{1}{2}e^{y/2}
\eeaa
as $n\to\infty.$ Combining the above three identities, we see that
\beaa
h_n\to \frac{1}{\sqrt{8\pi}}e^{y/2}
\eeaa
as $n\to\infty.$ Therefore, we conclude from Lemma \ref{childish} and (\ref{China}) that
\beaa
 \lim_{n\to\infty} P((n-2)\log (1-L_n^2) + 4\log p -\log \log p \geq y) = e^{-Ke^{y/2}}
\eeaa
for any $y\in\mathbb{R},$ where $K=\frac{1}{\sqrt{8\pi}}.$ Since $\varphi(y):=e^{-Ke^{y/2}}$ is continuous for all $y\in \mathbb{R},$ it is trivial to check that
\bea\lbl{post}
 \lim_{n\to\infty} P((n-2)\log (1-L_n^2) + 4\log p -\log \log p \leq y)=1- e^{-Ke^{y/2}}
\eea
for any $y\in\mathbb{R}.$ We get (\ref{watch}).\ \ \ \ \ \ \ \ $\blacksquare$\\

\noindent\textbf{Proof of Corollary \ref{sister}}. Dividing (\ref{river}) by $\log p$, we see that
\bea\lbl{foot}
 \frac{n}{\log p}\log (1-L_n^2) \to -4
\eea
in probability as $n\to\infty.$ By (i) of Theorem \ref{kind},  we know  $L_n\to 0$ in probability as $n\to\infty.$ Since $\rho_{ij}$ has density $f(\rho)$ as in (\ref{meeting}) for $i\ne j$, we have $P(L_n=0)=0$ for all $n\geq 3.$ Notice the function
\beaa
h(x):=\begin{cases}
x^{-1}\log (1-x), & \text{if $x\in (0,1)$;}\\
-1, & \text{if $x=0$}
\end{cases}
\eeaa
 is continuous on $[0,1),$ we have
\beaa
\frac{\log (1-L_n^2)}{L_n^2}=h(L_n^2) \to h(0)=-1
\eeaa
in probability as $n\to\infty.$ This together with  (\ref{foot}) yields
\beaa
\frac{n}{\log p}\cdot L_n^2 \to 4
\eeaa
in probability as $n\to\infty.$ The desired conclusion then follows.\ \ \ \ \ \ \ $\blacksquare$\\

\noindent\textbf{Proof of Corollary \ref{brother}}. By Theorem \ref{kind},
\bea\lbl{sugar}
P\left(n\log (1-L_n^2) + 4\log p -\log \log p \leq y\right) \to F(y)
\eea
as $n\to\infty,$ where $F(y)=1- e^{-Ke^{y/2}},\  y\in\mathbb{R}.$ Set
\bea\lbl{circle}
y_{n,p}=n\Big[1-\exp\Big\{\frac{1}{n}(-4\log p +\log \log p +y)\Big\}\Big].
\eea
Then, (\ref{sugar}) becomes that $P(nL_n^2 \geq y_{n,p})\to F(y),$ and hence
\bea\lbl{boy}
P(nL_n^2 -4\log p +\log\log p< y_{n,p}-4\log p +\log\log p)\to 1-F(y)
\eea
as $n\to\infty$ for any $y\in\mathbb{R}.$ We claim
\bea\lbl{plastic}
y_{n,p}-4\log p +\log\log p\to -(y + 8\alpha^2) \ \ \mbox{if}\ \ \frac{\log p}{\sqrt{n}}\to \alpha \in [0,\infty).
\eea
If this is true, by (\ref{boy}) and the continuity of $F(y),$
\beaa
\lim_{n\to\infty}P(nL_n^2 -4\log p +\log\log p \leq -(y + 8\alpha^2) )= 1-F(y)
\eeaa
for any $y\in \mathbb{R}.$ In other words, $nL_n^2 -4\log p +\log\log p$ converges weakly to a probability distribution function
\beaa
G(z):=1-F(-z-8\alpha^2)=\exp\{-Ke^{-(z+8\alpha^2)/2}\},\ z\in \mathbb{R},
\eeaa
as $n\to\infty.$ Now we  prove claim (\ref{plastic}).

In fact, set $t=-4\log p +\log \log p +y.$ Then $t=O(\log p)$ and $\frac{t}{n} \to 0 $ as $n\to\infty$ under the assumption $\frac{\log p}{\sqrt{n}}\to \alpha$. Consequently, by (\ref{circle}) and the Taylor expansion,
\beaa
y_{n,p}=n(1-e^{t/n})
& = & -n\Big[\frac{t}{n} +\frac{t^2}{2n^2} + O\Big(\frac{t^3}{n^3}\Big)\Big]\\
&=& -t-\frac{t^2}{2n}+ O\Big(\frac{t^3}{n^2}\Big)
\eeaa
as $n\to\infty.$ If  $\frac{\log p}{\sqrt{n}}\to \alpha$ as $n\to\infty,$ then $\frac{t^2}{2n}\to 8\alpha^2$ and $\frac{t^3}{n^2}\to 0$ as $n\to\infty.$ Therefore, (\ref{plastic}) is concluded. \ \ \ \ \ \ \ \ $\blacksquare$\\


\noindent\textbf{Proof of Theorem \ref{jean}}. (i). Assume (ii) of the theorem holds. Since $(\log p)/n \to \beta$ as $n\to\infty,$  dividing (\ref{divide}) by $n$, we see that $\log (1- L_n^2) \to -4\beta$ in probability, or equivalently, $L_n\to \sqrt{1-e^{-4\beta}}$ in probability as $n\to\infty.$

(ii). The proof here does not rely on the conclusion in (i). We first show that
\bea\lbl{eight}
(n-2)\log (1-L_n^2) + 4\log p -\log\log p
\eea
 converges weakly to the distribution function $F(y)=1- e^{-K(\beta) e^{y/2}},\  y\in\mathbb{R},$ where $K(\beta)$ is as in (\ref{east}). If this is true, by the condition $(\log p)/n \to \beta$ and the argument as in (i) above, we see that
\beaa
\log (1-L_n^2) \to -4\beta
\eeaa
in probability as $n\to\infty.$ Thus, by the Slusky lemma,
\beaa
& & n\log (1-L_n^2) + 4\log p -\log\log p\\
&= & \left[(n-2)\log (1-L_n^2) + 4\log p -\log\log p\right] + 2\log (1-L_n^2)
\eeaa
converges weakly to the distribution function $F(y)=1- e^{-K(\beta) e^{(y+8\beta)/2}},\  y\in\mathbb{R}.$ We now prove (\ref{eight}).

Fix $y\in \mathbb{R}.$ Let $N=n-2$ and $t=t_{n}\in [0,1)$ such that
\beaa
t^2=1-\exp\Big\{\frac{1}{N}(-4\log p +\log\log p +y)\wedge 0\Big\} .
\eeaa
It is easy to see that
\bea\lbl{USA}
P((n-2)\log (1-L_n^2) + 4\log p -\log\log p \geq y)=P(L_n\leq t)
\eea
as $n$ is sufficiently large, and
\bea\lbl{lake}
\lim_{n\to\infty}t_n= \sqrt{1- e^{-4\beta}}\in (0,1)\ \ \mbox{and}\ \ \ N\log (1-t^2)= -4\log p +\log\log p +y
\eea
as $n$ is sufficiently large.  We now calculate $h_n$ in Lemma \ref{childish} to obtain $\lim_{n\to\infty}P(L_n\leq t)$. Review
\bea\lbl{tame}
h_n = \frac{n^{1/2}p^2}{\sqrt{2\pi}}
\int_{t}^1(1-x^2)^{\frac{n-4}{2}}\,dx.
\eea
 It follows from Lemma \ref{also} and the first identity in (\ref{lake}) that
\beaa
n^{1/2}p^2
\int_{t}^1(1-x^2)^{(n-4)/2}\,dx &\sim & \frac{n^{1/2}p^2}{nt}\left(1-t^2\right)^{(n-2)/2}\\
& \sim & \frac{1}{\sqrt{1- e^{-4\beta}}}\cdot\frac{p^2}{\sqrt{N}}\left(1-t^2\right)^{N/2}
\eeaa
as $n\to\infty.$ By using the second identity in (\ref{lake}), we see that
\beaa
\frac{p^2}{\sqrt{N}}\left(1-t^2\right)^{N/2} & = & \frac{p^2}{\sqrt{N}}\cdot \exp\Big\{\frac{-4\log p + \log\log p +y}{N}\cdot \frac{N}{2}\Big\}\\
& = & \frac{\sqrt{\log p}}{\sqrt{N}}\cdot e^{y/2} \to \sqrt{\beta}\, e^{y/2}
\eeaa
as $n\to\infty.$ Collect all the facts above to have
\beaa
\lim_{n\to\infty}h_n= K(\beta) e^{y/2}
\eeaa
where
\beaa
K(\beta)=\Big(\frac{\beta}{2\pi(1-e^{-4\beta})}\Big)^{1/2}.
\eeaa
 By (\ref{USA}) then Lemma \ref{childish} we have
\beaa
 \lim_{n\to\infty} P((n-2)\log (1-L_n^2) + 4\log p -\log \log p \geq y) = e^{-K(\beta) e^{y/2}}
\eeaa
for any $y\in\mathbb{R}.$ By the same argument as getting (\ref{post}), the above yileds that
\beaa
 \lim_{n\to\infty} P((n-2)\log (1-L_n^2) + 4\log p -\log \log p \leq y)=1- e^{-K(\beta)e^{y/2}}
\eeaa
for any $y\in\mathbb{R}.$ We eventually arrive at (\ref{eight}).\ \ \ \ \ \ \ \ $\blacksquare$\\

\noindent\textbf{Proof of Theorem \ref{learn}}. (i). Assuming (ii) of the theorem, dividing (\ref{eve}) by  $\log p$, we see that
\beaa
\frac{n}{\log p}\log (1-L_n^2) \to -4
\eeaa
in probability as $n\to\infty.$ Since $(\log p)/n \to \infty$, we have $L_n\to 1$ in probability as $n\to\infty.$

(ii). The proof in this part does not rely on the conclusion in (i). Fix $y\in \mathbb{R}.$ Let $N=n-2$ and $t=t_{n}\geq 0$ such that
\beaa
t^2=1-\exp\Big\{\frac{1}{N}(-4\log p +\log n +y)\wedge 0\Big\}.
\eeaa
Obviously, $t_n \to 1^-$ as $n\to\infty$ by the condition $(\log p)/n \to \infty.$ Thus, without loss of generality, assume $t=t_n\in (0,1)$ for all $n\geq 1.$ Easily,
\bea
& & \log (1-t^2) =\frac{-4\log p +\log n +y }{N}\ \ \mbox{and}\lbl{mini}\\
& & P((n-2)\log (1-L_n^2) + 4\log p -\log n \geq y)=P(L_n\leq t)  \lbl{Russia}
\eea
as $n$ is sufficiently large. We now evaluate $h_n$ in Lemma \ref{childish} to obtain $\lim_{n\to\infty}P(L_n\leq t).$ Recall
\bea\lbl{user1}
h_n = \frac{n^{1/2}p^2}{\sqrt{2\pi}}
\int_{t}^1(1-x^2)^{\frac{n-4}{2}}\,dx.
\eea
From Lemma \ref{also} and the fact $t_n\to 1$ as $n\to\infty$ we obtain
\beaa
n^{1/2}p^2
\int_{t}^1(1-x^2)^{(n-4)/2}\,dx &\sim & \frac{n^{1/2}p^2}{nt}\left(1-t^2\right)^{(n-2)/2}\\
& \sim & \frac{p^2}{\sqrt{N}}\left(1-t^2\right)^{N/2}
\eeaa
as $n\to\infty.$ Combine this and (\ref{mini}) to have
\beaa
\frac{p^2}{\sqrt{N}}\left(1-t^2\right)^{N/2} & \sim & \frac{p^2}{\sqrt{N}}\cdot \exp\Big\{\frac{-4\log p + \log n +y}{N}\cdot \frac{N}{2}\Big\}\\
& = & e^{y/2}\cdot \frac{\sqrt{n}}{\sqrt{N}} \to e^{y/2}
\eeaa
as $n\to\infty.$ Joining all the above we have that
\beaa
\lim_{n\to\infty}h_n = \frac{1}{\sqrt{2\pi}}e^{y/2}
\eeaa
as $n\to\infty.$ From (\ref{Russia}) then Lemma \ref{childish} we finally obtain
\beaa
 \lim_{n\to\infty} P((n-2)\log (1-L_n^2) + 4\log p -\log n \geq y) = e^{-Ke^{y/2}}
\eeaa
for any $y\in\mathbb{R},$ where $K=\frac{1}{\sqrt{2\pi}}.$ By the same argument as getting (\ref{post}), the above actually implies  that
\beaa
 \lim_{n\to\infty} P((n-2)\log (1-L_n^2) + 4\log p -\log n \leq y)=1- e^{-Ke^{y/2}}
\eeaa
for any $y\in\mathbb{R}.$ This says that
\bea\lbl{liquor}
(n-2)T_n + 4\log p -\log n \Longrightarrow F(y)
\eea
with $F(y)=1- e^{-Ke^{y/2}},\  y\in\mathbb{R}$ and $K=1/\sqrt{2\pi}.$ Further, multiplying the left hand side of (\ref{liquor}) by $\frac{2}{n-2}$ we obtain
\bea\lbl{romance}
2T_n+ \frac{8\log p}{n-2}\overset{P}{\to} 0
\eea
as $n\to\infty.$ Noticing $(n-2)T_n + 2T_n=nT_n.$ Adding up (\ref{liquor}) and (\ref{romance}), we conclude from the Slusky lemma that
\beaa
nT_n + 4\log p -\log n +\frac{8\log p}{n-2} = nT_n + \frac{4n}{n-2}\log p-\log n
\eeaa
converges weakly to the distribution function $F(y)=1- e^{-Ke^{y/2}},\  y\in\mathbb{R}$ with $K=1/\sqrt{2\pi}.$
\ \ \ \ \ \ \ \ $\blacksquare$

\subsection{Proofs for Results on $\tilde{L}_n$ in Section \ref{silent}}\lbl{elm}

The proofs of the results on  $\tilde{L}_n$ are analogous to those of the results on $L_n$. The essential difference is to apply (\ref{meeting'}) in place of (\ref{meeting}).  Keeping all other arguments, we then get the proofs of the results on  $\tilde{L}_n$ stated in Section \ref{silent}. We omit the details for reasons of space.



\begin{thebibliography}{99}\lbl{referenceb}


\bibitem{Ahlfors}
Ahlfors, L. V. (1979). {\it Complex Analysis}, Third Edition. McGraw-Hill, New York.

\bibitem{AGG89}
Arratia, R. and Goldstein, L. and Gordon, L. (1989).
Two moments suffice for {P}oisson approximation: The Chen-Stein method.
{\it Ann. Probab. \bf 17}, 9-25.


\bibitem{BYT}
Bickel, P. J., Ritov, Y. and Tsybakov, A. B. (2009). Simultaneous analysis of Lasso and Dantzig selector.
{\it Ann. Statist. \bf  37}, 1705-1732.

\bibitem{CaiJiang}
Cai, T. T. and Jiang, T. (2010). Limiting Laws of Coherence of Random Matrices with Applications to Testing Covariance Structure and Construction of Compressed Sensing Matrices.
{\it Ann. Statist.}, to appear.

\bibitem{caiL}
Cai, T. T. and Lv, J. (2007). Discussion: The Dantzig selector:
statistical estimation when $p$ is much larger than $n$.
{\it Ann. Statist. \bf 35}, 2365-2369.

\bibitem{CaWaXu1}
Cai, T. T., Wang, L. and Xu, G. (2010a).
Shifting inequality and recovery of sparse signals,
{\it IEEE Transactions on Signal Processing \bf 58}, 1300-1308.

\bibitem{CaWaXu2}
Cai, T. T.,  Wang, L. and Xu, G.  (2010b).
Stable recovery of sparse signals and an oracle inequality,
{\it IEEE Transactions on  Information Theory \bf 56}, 3516-3522.

\bibitem{CanPlan}
Cand\`es, E. J. and  Plan, Y. (2009).
Near-ideal model selection by $\ell_1$ minimization.
{\it Ann. Statist. \bf 37},  2145-2177.

\bibitem{CanTao07}
Cand\`es, E. J. and  Tao, T. (2007). The Dantzig selector: statistical
estimation when $p$ is much larger than $n$ (with discussion).
{\it Ann. Statist. \bf 35}, 2313-2351.

\bibitem{CanTao05}
Cand\`es, E. J. and  Tao, T. (2005). Decoding by linear programming,
{\it IEEE Trans. Inf. Theory \bf 51}, 4203-4215.

\bibitem{Donoho06}
Donoho, D. (2006). Compressed sensing.
{\it IEEE Trans. Inf. Theory \bf 52}, 1289-1306.

\bibitem{}
Donoho, D. L., Elad, M. and Temlyakov, V. N. (2006). Stable recovery
of sparse overcomplete representations in the presence of noise.
{\it IEEE Trans. Inf. Theory \bf 52}, 6-18.

\bibitem{DonHuo}  Donoho, D. L. and Huo, X. (2001).
Uncertainty principles and ideal atomic decomposition.
{\it IEEE Trans. Inf. Theory \bf 47}, 2845-2862.

\bibitem{FanL08}
Fan, J. and Lv, J. (2008). Sure independence screening for ultrahigh
dimensional feature space (with discussion).
{\it J. Roy. Statist. Soc. Ser. B} 70, 849-911.

\bibitem{FL10}
Fan, J. and Lv, J. (2010).
A selective overview of variable selection in high dimensional feature
space.
{\it Statistica Sinica} 20, 101-148.

\bibitem{Fuchs1}
Fuchs, J. J. (2004).
On sparse representations in arbitrary redundant bases.
{\it IEEE Trans. Inf. Theory \bf 50}, 1341-1344.

\bibitem{Gamelin}
Gamelin, T. W. (2001). {\it Complex Analysis.} Springer.

\bibitem{Jiang04}
Jiang, T. (2004). The asymptotic distributions of the largest entries of sample correlation
matrices. {\it Ann. Appl. Probab. \bf 14}, 865-880.

\bibitem{Eaton}
Kariya, T. and Eaton, M. L. (1977). Robust tests for spherical symmetry.
{\it Ann. Statist. \bf 5}, 206-215.

\bibitem{LLR08}
Li, D., Liu, W. and Rosalsky, A. (2009). Necessary and sufficient conditions for the asymptotic distribution of the largest entry of a sample correlation matrix.
{\it Probab. Theory Relat. Fields  \bf 148}, 5-35.

\bibitem{Li2010}
Li, D., Qi, Y. and Rosalsky, A. (2010).  On Jiang's asymptotic distribution of the largest entry of a sample correlation matrix. Preprint.

\bibitem{LR}
Li, D. and Rosalsky, A. (2006). Some strong limit theorems for the largest entries of sample correlation matrices.
{\it Ann. Appl. Probab. \bf 16}, 423-447.

\bibitem{LLS}
Liu, W., Lin, Z. and Shao, Q. (2008). The asymptotic distribution and Berry--Esseen bound of a new test for independence in high dimension with an application to stochastic optimization.
{\it Ann. Appl. Probab. \bf 18}, 2337-2366.


\bibitem{Muirhead1982}
Muirhead, R. J. (1982). {\it Aspects of Multivariate Statistical Theory.}
Wiley, New York.

\bibitem{Wigner58}
Wigner, E. P. (1958). On the distribution of roots of certain symmetric matrices.
{\it Ann. Math. \bf 67}, 325-328.

\bibitem{Zhou}
Zhou, W. (2007). Asymptotic distribution of the largest off-diagonal entry of correlation matrices.
{\it Transaction of American Mathematical Society \bf 359}, 5345-5363.
\end{thebibliography}
\end{document}